\documentclass[letterpaper,conference,10pt]{IWMRRF} \pdfoutput=1
\usepackage[pdftex]{graphicx}
\usepackage{pgf}
\usepackage{url}      
\usepackage{amsmath}                                           
\usepackage{amssymb}
\usepackage{color}
\usepackage{colortbl}             
                        
\usepackage{array}

\renewcommand{\d}{\mathrm{d}}
\begin{document}

\title{Stretching-Based Diagnostics in a Differential Geometry Setting}
\author{\authorblockN{Johannes Poppe\authorrefmark{1},
Dirk Lebiedz\authorrefmark{1}}
\authorblockA{\authorrefmark{1}Institute for Numerical Mathematics, Ulm University, Germany}}
\maketitle

\begin{abstract}
It is a prominent challenge to analytically characterize slow invariant manifolds for dynamical systems with multiple time-scales.
To this end, we transfer the system into a differential-geometric framework. %to identify these objects by local properties of the dynamics.
This setting enables to formulate stretching-based diagnostics in a new context, coinciding with the intrinsic differential-geometric property of sectional curvature.
\end{abstract}

\section{Introduction}
The identification of slow invariant manifolds (SIMs) is an essential part in model-order reduction for reactive systems. 
The mathematical definition of the SIM by Fenichel can be considered unsatisfactory, because it is only applicable to so-called slow-fast system and does not provide the uniqueness of the SIM. Observing the phase space of the dynamical system (not necessarily a slow-fast system), the SIM becomes a geometric object which attracts trajectories, resulting in a bundling behavior. We aim to find a more general definition of the SIM, guided by the prior observations in phase space within the field of differential geometry. This setting provides one major benefit: All quantities are formulated covariantly, i.e. they are independent of the coordinate choice.
A recent work by Heiter and Lebiedz \cite{heiter} translates the invariance property to vanishing sectional curvatures in the extended phase space. 
\section{Geodesics in Spacetime} 
Let $\dot{x} = f(x)$ with $x\in \mathbb{R}^n$ and $f:\mathbb{R}^n \to \mathbb{R}^n$ be sufficiently smooth. We consider the extended state space, i.e., the original state space is extended by an additional time-axis $\tau$. In order to clarify that $\tau$ is part of the extended state space, we call $\tau$ "explicit time".
In contrast, $t$ is called implicit time. The resulting curves in the extended system are the solutions of
\begin{equation}\label{eq:ext_sys}
\frac{\d}{\d t}
 \begin{pmatrix}
x \\
\tau
\end{pmatrix} = \begin{pmatrix}
f(x)\\1
\end{pmatrix} \in \mathbb{R}^{n+1}.
\end{equation}
Bundling behavior of trajectories of the original system corresponds to bundling of the solutions of $\eqref{eq:ext_sys}$. 
The space $M:=\mathbb{R}^{n+1}$ is trivially a manifold. By defining a metric tensor $g$ (a family of inner products that varies smoothly from point to point), the tuple $(M,g)$ becomes a Riemannian manifold. The core idea is to couple the dynamics of \eqref{eq:ext_sys} with a metric such that every trajectory becomes a geodesic - a shortest connection path for each tuple of points on the trajectory with regard to the chosen metric. An evident choice for the connection is the so-called Levi-Civita connection. In this setting, solutions of the extended system can be interpreted in analogy to free falling particles in a graviational field within the framework of general relativity. In the context of chemical reaction mechanisms, the gravitational field correlates to an abstract chemical force. A suitable choice for a metric $g$ can be derived and expressed by its components regarding the standard coordinates $x_1,\cdots,x_n,\tau$:
\[
g_{ij} = \begin{pmatrix}
\text{Id}_n  & -f(x) \\
-f(x)^T &  1+f(x)^T f(x) \\
\end{pmatrix}.
\]
The bundling trajectories of the extended system \eqref{eq:ext_sys} becomes a set of geodesics which bundle alongside a specific subset/submanifold of geodesics - the SIM in space time. Hence, the SIM is supposed to be characterized by some differential geometric property representing this bundling behavior.
\section{Geodesic Stretching}
In general relativity, bundling behavior is related to geodesic deviation, neighboring geodesics experience relative accelerations towards each other.
\subsection{Geodesic Deviation}\label{sec:geod_dev}
For $p \in M$, each geodesic $\ell : (-\varepsilon, \varepsilon) \to M$ passing through $p$ at $t = 0$ defines the tangent vector
\[
T :=\frac{\d \ell}{\d t} (0) \in T_p M, \qquad \text{here:}\qquad T=\begin{pmatrix}f(x)\\1,
\end{pmatrix}
\]
because we plug-in solution trajectories of $\eqref{eq:ext_sys}$. 
\begin{figure}[h]
	\centering
	\def\svgwidth{0.48\textwidth}
	\fbox{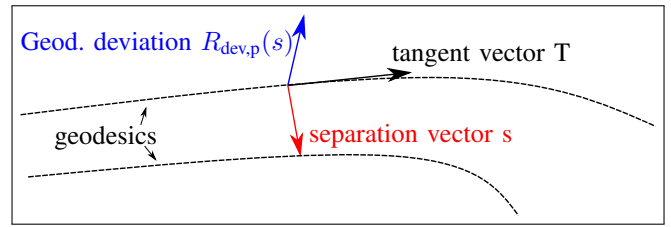}
	\caption{Visualization of geodesic deviation} \label{fig:geo_dev}
\end{figure}
Geodesic deviation is formally defined as the endomorphism $R_{\text{dev},p}:  T_p M \to T_p M$ formed by inserting the tangent vector $T$ into the first and third argument of the Riemann curvature tensor $R_p:(T_p M)^3 \to T_p M$,
\begin{equation} \label{eq:geod_dev}
R_{\text{dev},p}(s) := R_p(T,s)T.
\end{equation}
The intuitive concept of geodesic deviation is shown in Fig.$~\ref{fig:geo_dev}$. The input vector $s$ is the so-called "separation vector" representing a small difference between two points on the paths of two neighboring geodesics.
 The output vector represents the relative acceleration between both geodesics.
\subsection{Stretching-Based Diagnostics}
In \cite{adrover}, Adrover et al. study SIMs by comparing so-called stretching rates $\omega_v(x)$ of solution trajectories of a given dynamic system $\dot{x}=F(x)$ for different vectors $v(x) \in \mathbb{R}^n = T_x \mathbb{R}^n$. Each stretching rate is defined by
\[
\omega_v(x) := \frac{\langle J_Fv , v \rangle}{\langle v , v \rangle},
\]
where $J_F$ represents the Jacobian matrix of $F$. For a given sub-manifold $U \subset \mathbb{R}^n$ and $p \in U$, the tangent space $T_p \mathbb{R}^n$ is decomposed into the direct sum of tangent space $T_p U$ and normal space $N_p U$. % $T_p \mathbb{R}^n = T_p U \oplus N_p U.$
Adrover et al. argue that on the SIM, stretching in normal directions is supposed to dominate tangent directions. Hence - in case of a one-dimensional SIM in a two-dimensional system, with $v_t$ and $v_o$ being tangent and normal to the trajectory respectively - the so-called stretching ratio
\[
r(x) := \frac{\omega_{v_o}(x)}{\omega_{v_t}(x)}
\]
is supposed to be larger than one at the SIM. Bundling of trajectories is directly correlated to the stretching ratio. It appears evident that $r(x)$ is particularly large on the SIM, making its maximization on the phase space for a fixed choice of reaction progress variables (RPVs) a viable option to approximate the SIM.
\subsection{Adapting the Stretching-Based Approach}
We integrate this stretching-based approach within the framework from subsection \ref{sec:geod_dev} in the following way: We replace the euclidean metric $\langle \cdot, \cdot \rangle$ with the metric $g$ defined above and expand each $v\in T_p \mathbb{R}^n$ by an additional explicit time component (which is set to 0):
\[
 \mathbb{R}^n \ni v \mapsto \begin{pmatrix}
 v\\0
 \end{pmatrix}:= \tilde{v} \in \mathbb{R}^{n+1}.
\]
Instead of the Jacobian matrix acting on each tangent vector representing an infinitesimal acceleration, we apply the geodesic deviation endomorphism defined in \eqref{eq:geod_dev}. The result is the notion of so-called geodesic stretching
\[
\omega_{g,\tilde{v}}:= \frac{g(R_{\text{dev},p}(\tilde{v}),\tilde{v})}{g(\tilde{v},\tilde{v})}.
\]
\subsection{Geometrical Interpretation}
One can show that $\omega_{g,v}$ has an intrinsic differential geometric meaning: It represents the sectional curvature $K(\tilde{v}, w)$ of the expanded vector $\tilde{v}$ with tangent vector $T$ of each space time trajectory from above: 
\[
\frac{g(R_{\text{dev},p}(\tilde{v}),\tilde{v})}{g(\tilde{v},\tilde{v})} = \frac{g(R(T,\tilde{v})T,\tilde{v})}{g(\tilde{v},\tilde{v})g(T,T)-g(\tilde{v},T)^2} := K(\tilde{v},T).
\]
Sectional curvature depends on the given metric $g$. Hence, there is no direct connection to the utilized sectional curvature in \cite{heiter} which is based on a different metric.
\subsection{Results}
\begin{figure}[h]
\centering
\def\svgwidth{\columnwidth}
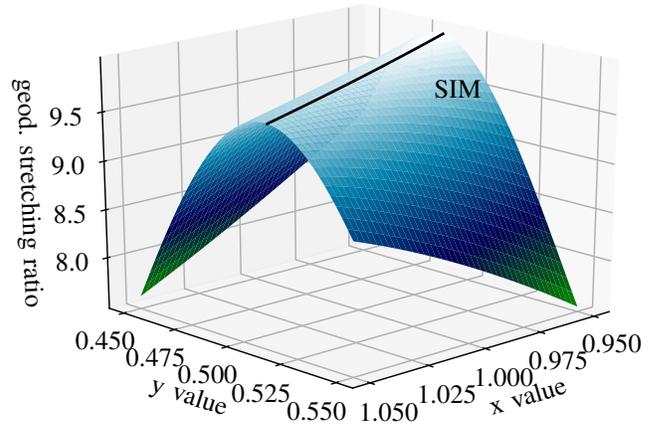
\caption{Surface plot of geodesic stretching ratio for the Davis-Skodje system near the SIM. Black curve is stretching ratio on SIM.}
\label{fig:stret_ratio_plot}
\end{figure}

We test our ansatz by considering the well-known non-linear Davis-Skodje (DS) system:
%\begin{align}
\begin{subequations}\label{eq:ds}
	\begin{align} 
		\frac{\d x}{\d t} &= -x \\
		\frac{\d y}{\d t} &= -\gamma y + \frac{(\gamma-1)x + \gamma x^2}{(1+x)^2}
	\end{align}
\end{subequations}
where $\gamma > 1$ is a fixed parameter and $x > 0$. This system has a one-dimensional SIM given by the graph representation
$
y = h(x) = \frac{x}{1+x}.
$
Fig. \ref{fig:stret_ratio_plot} shows geodesic stretching ratios near the SIM for system \eqref{eq:ds} with $\gamma=3$. Setting $x$ as RPV and maximizing the ratio with respect to $y$ yields points near the SIM.
\section*{Acknowledgement}
Special thanks goes to the Klaus-Tschira foundation for financial funding of the project.


\begin{thebibliography}{10}
	
\bibitem{heiter} P. Heiter and D. Lebiedz,  ``Towards Differential Geometric Characterization of Slow Invariant Manifolds
in Extended Phase Space: Sectional Curvature and Flow Invariance'', in {\em SIAM J. Applied Dynamical Systems}
Vol. 17 (2018), No. 1, pp. 732--753

\bibitem{adrover} A. Adrover, F. Creta, M. Giona, M. Valorani, ``Stretching-based diagnostics and reduction of chemical kinetic models with diffusion'', in {\em Journal of Computational Physics} 225 (2007) 1442--1471

\end{thebibliography}
\end{document}